\documentclass[leqno,12pt]{article}
\begin{document}

\noindent
{\large \bf On the vanishing of higher syzygies of curves. II}
\\\\
{\bf Marian Aprodu}\footnote{Supported by an E.C. Marie Curie Fellowship,
contract number HPMF-CT-2000-00895}
\\\\
{\small Romanian Academy, Institute of Mathematics "Simion Stoilow",
P.O.Box 1-764, RO-70700, Bucharest, Romania (e-mail:
Marian.Aprodu\char64 imar.ro) \&
\\
Universit\'e de Grenoble 1,
Laboratoire de Math\'ematiques, 
Institut Fourier BP 74,
38402 Saint Martin d'H\`eres Cedex,
France (e-mail: aprodu\char64 mozart.ujf-grenoble.fr)}

\bigskip

\bigskip

\noindent
One of the remarkable results of Segre's, quoted in \cite{AC} as Theorem
0.2, states that generic $k$-gonal curves have distinguished nodal models lying
on the Hirzebruch surface $\Sigma_1$, in such a way that minimal pencils are
given by the ruling. Since there exist several results relating Koszul
cohomology of a surface to the Koszul cohomology of curves which lie on
it, it seems natural to make use of this very geometric context to
verify some cases of Green's conjecture. 
As exemplified by the main result of this short Note (to be compared to
\cite{Sch}), this strategy works well for curves of gonality much smaller than
the generic value (as for the case of generic curves of large gonality, other
recently developped approaches, like the innovative one of \cite{Vo}, are much
more efficient).
\\\\
{\bf Theorem 1.}{\em 
Green's conjecture is valid for a generic $k$-gonal curve of genus
$g$, for which $g\geq k(k-1)/2$.}
\\\\
The idea of proof consists of showing some vanishing result for the
Koszul cohomology of blowups of $\Sigma _1$, and bringing this into
relation with Koszul cohomology of strict transforms of nodal curves on
$\Sigma _1$. The fact of working over $\Sigma _1$ is not essential, and we
can actually state similar results for an arbitrary Hirzebruch surface
$\Sigma_e$, as follows (we freely use the notation of \cite{Ap}, such as
$C_0$ for the minimal section, $f$ for the fibre of the ruling etc,
without further mention).
\\\\
{\bf Theorem 2.}{\em 
Let $e\geq 0$, $\alpha\geq 2$, $\beta$, and $\gamma$ be four integers
which satisfy the following inequalities $\beta\geq \mbox{\rm max}
\{ \alpha e, \alpha+e \}$, and $\gamma\leq \beta-(e-1)\alpha$.
Let $\Sigma _\Gamma\stackrel{\sigma_\Gamma} {\rightarrow}\Sigma_e$
be the blowup of $\Sigma_e$ in a set $\Gamma=\{x_1,...,x_\gamma\}$ 
of generic points of $\Sigma_e$, denote $E_\Gamma$ the exceptional divisor,
and $H_\Gamma=\sigma_\Gamma^*{\cal O}_{\Sigma_e}(\alpha C_0+\beta f)-E_\Gamma$.
Then $K_{p,1}(\Sigma _\Gamma,H_\Gamma)=0$ for all $p\geq h^0
H_\Gamma -\alpha -1$.}
\\\\
{\em Proof.} We show first that {\em there exists} a set of points
$\Gamma=\{x_1,...,x_\gamma\}$, for which vanishing of the Koszul
cohomology on $\Sigma _\Gamma$ holds as stated in the conclusion.
We denote $L={\cal O}_{\Sigma_e}(\alpha C_0+\beta f)$, 
and we choose $Y\in |L|$ a smooth curve (such a curve does exist, as
$\beta \geq\alpha e$, and $\beta>0$). Since $\gamma\leq
\beta-(e-1)\alpha=\mbox{deg}\big({\cal O}_{\Sigma_e}(C_0+f)_{|Y}\big)$,
there exists a set of points $\Gamma=\{x_1,...,x_\gamma\}
\subset Y$, such that ${\cal O}_{\Sigma_e}(C_0+f)_{|Y}-x_1-...-x_\gamma$
is effective.
Therefore, Theorem 6.3 of \cite{Ap} applies in our case, 
and thus   $K_{p,1}(Y,L_{|Y}-x_1-...-x_\gamma)=0$
for all $p\geq h^0(Y,L_{|Y}-x_1-...-x_\gamma)-\alpha$.
In particular, by means of \cite{Ap} Remark 1.3 applied
for the strict transform $\widetilde{Y}\in |H_\Gamma |$ of $Y$, we obtain 
$
K_{p,1}(\Sigma_\Gamma,H_\Gamma)=0,
$ for all 
$
p\geq h^0(\Sigma_\Gamma,H_\Gamma)- \alpha -1
$, 
as claimed. Furthermore, as the line bundle $L_{|Y}-x_1-...-x_\gamma$ is
nonspecial, the points $\{x_1,...,x_\gamma\}$ impose independent
conditions on the linear system $|L|$.
\\
Next, we claim the vanishing of the Koszul cohomology of bundles
$H_\Gamma$ on different blowups of $\Sigma_e$ is an open condition in 
families of $\Gamma$'s which impose independent conditions on $|L|$. This
comes from the fact that we can describe
$K_{p,1}(\Sigma_\Gamma,H_\Gamma)$ as the cohomology in the middle of
the complex 
$$
\bigwedge^{p+1}H^0(H_\Gamma)
\longrightarrow
H^0(H_\Gamma)\otimes\bigwedge^pH^0(H_\Gamma) 
\longrightarrow
H^0(L^{\otimes 2})\otimes\bigwedge^{p-1}H^0(H_\Gamma) .
$$
As the first map is always injective, and the complex above is the 
fiber-restriction of a complex of vector bundles over the locus
parametrizing the cycles $\Gamma$ which impose independent conditions
on $|L|$, the claim follows.
\\\\
{\em Remark 3. } The condition $\beta\geq \mbox{max}
\{ \alpha e, \alpha+e \}$ is equivalent to saying that the ruling
on $\Sigma_e$ restricts to a minimal pencil to any irreducible smooth
curve in the linear system $|L|$.
\\\\
An immediate consequence of Theorem 2 is the following.
\\\\
{\bf Corollary 4.}{\em 
Let $e\geq 0$, $k\geq 4$, $m$, and $\gamma$ be four integers which satisfy 
the following inequalities
$m\geq \mbox{\rm max}\{ (k-1)e+2, k+2e \}$,
and $\gamma\leq m-e-2-(k-2)(e-1)$. Let $\Gamma=\{x_1,...,x_\gamma\}$ be a 
set of generic points of $\Sigma_e$, and let $X$ be an irreducible  curve on
$\Sigma_e$, numerically equivalent to $kC_0+mf$, having ordinary nodes at
$x_1,...,x_\gamma$ and no other singular points. If $\widetilde{X}$ denotes 
the normalization of $X$, then the Clifford dimension of $\widetilde{X}$ 
equals one, and Green's conjecture is valid for $\widetilde{X}$.}
\\\\
{\em Proof.} The proof runs in a similar way to that of \cite{Ap},
8.1, so we shall only sketch it here. We set $\alpha =k-2$, and 
$\beta = m-e-2$, and we use the same notation as in the proof of Theorem 2. 
We also denote the genus of $\widetilde{X}$ by
$g=(k-1)(m-1-ke/2)-\gamma$. We observe that $h^0(H_\Gamma-\widetilde{X})=
h^1(H_\Gamma-\widetilde{X})=0$, and $H_{\Gamma|\widetilde{X}}=
K_{\widetilde{X}}$, which altogether yield to a long exact sequence:
$$
...
\longrightarrow
K_{p,1}(\Sigma_\Gamma,H_\Gamma)
\longrightarrow 
K_{p,1}(\widetilde{X},K_{\widetilde{X}})
\longrightarrow
K_{p-1,2}(\Sigma_\Gamma,-\widetilde{X},H_\Gamma)
\longrightarrow
...
$$
By means of Green's vanishing \cite{Gr} 3.a.1., we obtain 
$K_{p-1,2}(\Sigma_\Gamma,-\widetilde{X},H_\Gamma)=0$ for 
$p\geq h^0(2H_\Gamma -\widetilde{X})+1$. As 
$h^0(2H_\Gamma -\widetilde{X})+1\leq g-k+1$, we apply Theorem 2 to conclude.
\\\\
{\em Proof of Theorem 1.} 
Thanks to \cite{AC}, 4.7, we know that for any integers $k$, $m$, and
$\gamma$ satisfying the inequalities $k\geq 4$, $m\geq k+2$, and $\gamma\leq
m-3$, and for any set $\Gamma$ of $\gamma$ points in general position on
$\Sigma_1$, there exists an irreducible curve $X$ on $\Sigma_1$, numerically
equivalent to $kC_0+mf$, and having ordinary nodes at the points of $\Gamma$
and no other singular points. 
In view of Corollary 4, and of semicontinuity of graded Betti numbers,
what is left from the proof is now a purely numerical
matter: for any integers $k\geq 4$, and $g\geq k(k-1)/2$, there exist $m\geq
k+2$, and $0\leq\gamma \leq m-3$, such that $g=(k-1)(m-1-k/2)-\gamma$.

\end{document}